\documentclass{amsart}

\usepackage{amsmath, amsfonts,amssymb,graphicx,amsthm,color}

\newtheorem{theorem}{Theorem}
\newtheorem*{theorem*}{Theorem}

\newtheorem*{corollary*}{Corollary}

\newtheorem*{definition*}{Definition}
\newtheorem{lemma}{Lemma}
\newtheorem*{lemma*}{Lemma}

\newtheorem{proposition}{Proposition}
\newtheorem*{proposition*}{Proposition}

\numberwithin{equation}{section}
\newtheorem{remark}{Remark}
\newtheorem{conjecture}{Conjecture}

\def\R{{\mathbb R}}

\def\E{{\mathbb E}}
\def\Tr{{\text{Tr}}}
\def\P{{\mathbb P}}
\def\cala{{\mathcal{A}}}
\def\<{{\langle}}
\def\>{{\rangle}}
\def\Var{{\text{Var}}}

\title[The spectrum of 
 random kernel matrices]{The spectrum of random kernel matrices: universality results for rough and varying kernels}
\author{Yen Do}
\email{yen.do@yale.edu}
\thanks{Y. Do is partially supported by NSF grant DMS-1201456.}

\author{Van Vu}
\email{van.vu@yale.edu}
\address{Department of Mathematics, Yale University, New Haven, CT 06520-8283}
\thanks{V. Vu is partially supported by NSF grant  DMS-0901216}

\subjclass[2000]{15B52}

\date{\today}

\begin{document}

\begin{abstract}
We consider random matrices whose entries are $f(X_i^TX_j)$ or $f(\|X_i-X_j\|^2)$ for  iid vectors $X_i \in \R^p$ with normalized distribution.  Assuming that  $f$ is sufficiently smooth and the distribution of $X_i$'s is sufficiently nice, El Karoui \cite{karoui} showed that the spectral distributions of these matrices behave  as if $f$ is \emph{linear} in the Mar\v{c}henko--Pastur limit. When $X_i$'s are Gaussian vectors, variants of this phenomenon were recently proved for varying kernels, i.e. when $f$ may depend on $p$, by Cheng--Singer \cite{cheng-singer}.  Two results are shown in this paper: first it is shown that for a large class of distributions the regularity assumptions on $f$ in El Karoui's results can be reduced to minimal; and secondly it is shown that the Gaussian assumptions in Cheng--Singer's result can be removed, answering a question posed in \cite{cheng-singer} about the universality of the limiting spectral distribution. 
\end{abstract}

\maketitle

\section{Introduction}
Let  $X_1, \dots, X_n \in \mathbb R^p$ be iid random vectors with normalization $\E [X_i] = 0$ and $\E \|X_i\|^2 = 1$, here $\|.\|$ denotes the Euclidean length in $\R^p$. For any $F: \R^p \times \R^p \times \R \rightarrow \R$ symmetric in the first two variables, consider the random matrix $A$ with entries
\begin{equation}\label{e.Adefdiag}
A_{ij} = F(X_i , X_j,p) \ \ ,
\end{equation}
or the variant with zeros on the diagonal
\begin{equation}\label{e.Adefnodiag}
\qquad A_{ij} = \begin{cases} F(X_i , X_j,p) , & i\ne j \\
0, & i=j \ .
\end{cases}
\end{equation}
Following previous literature \cite{karoui,cheng-singer}, in this paper these matrices will be refered to as \emph{random kernel matrices} generated by $F$ and the distribution of $X_i$'s. As described in \cite{karoui}, practical examples of $F$ are of the form
$$F(X,Y,p) = f(X^T Y, p) \ \ \text{or} \ \ f(\|X-Y\|^2,p) \ \ ,$$
More generally, one could have $F(X,Y,p)=f(g(X,Y),p)$ for some symmetric $g:\R^p \times \R^p \to \R$; some Lemmas in this paper are stated in this setting under suitable normalizing assumptions on $g$ (relative to $X_i$'s). For convenience, $g$ will be refered to as the  \emph{kernel} and $f$ will be refered to as the \emph{envelope} that generate $A$.  Examples of envelope functions are $f(x) = \exp(xa),  f(x)=(1+x)^a$, where $a$ is fixed; others can be found in Rasmussen--Williams \cite{rasmussen-williams} and Williams--Seeger \cite{williamsseeger}. 

We will be interested in weak-limit of the empirical distribution $\rho_A$ of $A$ when $p,n\to \infty$ such that $p/n \to  \gamma \in (0,\infty)$, a fixed constant. Recall that
$$\rho_A(x):= n^{-1}\sum_{i=1}^n \delta_{\lambda_i}(x) \ \ ,$$
here $\lambda_1,\dots, \lambda_n$ are eigenvalues of $A$ and $\delta _{\lambda}$ denotes the counting measure at $\lambda$. This research direction has been investigated recently by El Karoui \cite{karoui} and Cheng--Singer \cite{cheng-singer}, motivated by studies from machine learning and statistical analysis. In \cite{karoui}, it was assumed that $A$ is generated by either the inner-product or the distance kernels, and with $p$-independent envelope functions (which is the natural setting relative to the above normalization of $X_i$). It was shown in \cite{karoui} that for $f$ sufficiently smooth the limiting behavior of $\rho_A$ depends only on a {\it linear} component of $f$. It turns out that a variant of this phenomenon continues to hold even if $f$ depends on $p$: for $g(X,Y)=X^TY$ this was proved for Gaussian random vectors in a recent result of Cheng and Singer \cite{cheng-singer}. See also Bordenave \cite{bordenave12} for a related recent work in this direction that appeared after an initial circulation of a first draft of this paper.

The goal of this paper is to extend the results in \cite{karoui,cheng-singer} to more general settings. In particular, Theorem~\ref{t.main} will (positively) answer a question by Cheng and Singer \cite{cheng-singer}
regarding the universality of the limiting spectral distribution of Cheng--Singer's models. 

We would like to point out that El Karoui \cite{karouinoise} also considered a related model where the entries $g(X_i^TX_j)$ are affected by random noise before the envelope $f$ is applied outside; the reader is refered to the beautiful work \cite{karouinoise} for further details. There is also a vast amount of literature concerning limiting behaviors of $\rho_A$ when $p$ is low or fixed, the interested reader is refered to  \cite{mezardetal,vershik,williamsseeger,koltchinskiigine,bogomolnyetal, bordenave, jiang} and references there-in.

For clarity, below the discussion of previous and new results and are divided into two sections.

\subsection{The $p$-independent setting}
In this section, the setting when $F$ is independent of $p$ (relative to the above normalization of $X_i$) will be discussed. In other words, only the settings when $F(X,Y,p) = f(X^T Y)$ or $f(\|X-Y\|^2)$ (for some $p$-independent envelope function $f$) will be considered. Since the vectors $X_i$ are normalized, this is the natural setting for $f$.

\subsubsection{The inner product kernel}
Let $F(X,Y,p)=f(X^T Y)$. When the limiting spectral distribution for the model \eqref{e.Adefdiag} of $A$ is known, standard arguments may be used to deduce   the limiting spectral distribution of \eqref{e.Adefnodiag} (see e.g. \cite{bai-silverstein} or \cite{karoui}, see also Lemma~\ref{l.removediagonal} of the current paper), and vice versa. Below the model \eqref{e.Adefdiag} will be assumed unless otherwise stated.

For \emph{linear} envelope functions, it is well-known that if the distribution of the vectors $X_i$'s satisfies certain martingale/concentration properties then $\rho_A$ converges weakly to some form of the Mar{\v{c}}enko--Pastur (MP) distribution, whose density is given by
 $$\rho_{MP,\gamma}(x):= (1- \gamma)1_{\{\gamma <1\}} \delta_0(x) + \frac {\gamma}{2\pi x} \sqrt{(b-x)(x-a)} 1_{[a,b]}(x) \ \ ,$$
$$a:=(1-\frac 1{\sqrt \gamma})^2 \qquad \text{and} \qquad b:=(1+\frac1{\sqrt \gamma})^2 \ \ .$$
This convergence was first established by Mar{\v{c}}enko--Pastur \cite{marchenkopastur} (see also Wachter \cite{wachter}) when the entries of each vector $X_i$ are iid. Various authors have then extended this result to more general settings, see e.g.  Auburn \cite{auburn}, Yin and Krisnaiah \cite{yin-krisnaiah}, Silverstein \cite{silverstein95}, G\"otze and Tikhomirov \cite{gotze-tikhomirov-04, gotze-tikhomirov-05}, El Karoui \cite{karoui07a, karoui}, Adamczak \cite{adamczak}, Pajor and Pastur \cite{pajor-pastur}, Bordenave et~al. \cite{bordenaveetal}, Chafai \cite{chafai}, Chatterjee et~al. \cite{chatterjeeetal}, and O'Rourke \cite{orourke}.  In particular, the result holds for $X_i$'s drawn independently from isotrophic log-concave distributions, this is a result of Pajor and Pastur \cite{pajor-pastur}. Extensions to settings with some martingale-type assumptions were carried out in \cite{gotze-tikhomirov-04, gotze-tikhomirov-05, adamczak} , and extensions to settings with some concentration conditions on the distributions of $X_i$'s were done in \cite{karoui07a}. See also \cite{orourke, bordenaveetal, chatterjeeetal, chafai} for other generalizations.

For \emph{nonlinear} envelope functions with sufficient smoothness, it was shown by El Karoui \cite{karoui} that if the distribution of $X_i$'s is sufficiently nice then $A$ has the same the limiting spectral distribution as
$$B= [f(1)-f(0)-f'(0)]I_n + f'(0) (X_i^TX_j)_{i,j} \ \ .$$
Here and in the rest of the paper, $I_n$ is the identity $n\times n$ matrix. Let $\alpha:=f(1) - f(0)-f'(0)$. Using the linear theory, it follows that
\begin{equation}\label{e.karoui-inner}
\rho_A(x) \stackrel{weak}{\to}  \begin{cases}
\frac 1 {|f'(0)|} \rho_{MP,\gamma}(\frac {x-\alpha}{f'(0)}) , & \text{if $f'(0)\ne 0$;}\\
\delta_\alpha, & \text{if $f'(0)=0$.}
\end{cases}
\end{equation}
In El Karoui \cite{karoui}, the convergence \eqref{e.karoui-inner} was considered in two different settings:

(i) The iid setting with $K$ moment bounds: Assume that  
\begin{equation}\label{e.iid}
\begin{cases}
\text{the entries of $X_i$ are iid with}\\
\E |\sqrt p X_{ij}|^{K} = O(1) \ \ .
\end{cases}
\end{equation} 
In this setting,  it was shown in \cite{karoui} that \eqref{e.karoui-inner} holds if: $K>4$ and $f$ is $C^3$ near $0$ and $C^1$ near $1$.

(ii) The high concentration setting with parameter $c(p)$: Assume that 
\begin{equation}\label{e.karouiconcentration}
\begin{cases}
\text{for any $1$-Lipschitz function $F$}\\
\text{there exists absolute constants $C,b>0$ such that}\\
P(|F(X_i)-m_F|>t) \le C \exp(-c(p) t^b)   \ \ , \ \ \forall t>0  \ \ , \\
\text{here $m_F$ denotes a median of $F(X_i)$.}
\end{cases}
\end{equation}
In this setting, it was shown in \cite{karoui}  that \eqref{e.karoui-inner} holds if two conditions hold: 
\begin{itemize}
\item $f$ is $C^2$ near $0$ and   $C^1$ near $1$, and 
\item $c(p) \ge Cp^{\epsilon+b/4}$ for some absolute constant $C>0$. (For simplicity we'll write $c(p) \ge O(p^{\epsilon+b/4})$.)
\end{itemize}

In the iid setting \eqref{e.iid}, it was shown in \cite{karoui} that a stronger convergence in spectral norm  holds. In particular, it follows that some information about  the largest  eigenvalue  of $B$ could be transfered to $A$. The interested reader is refered to \cite{karoui, karoui07b} and the references there-in for related literature.

The estimate \eqref{e.karouiconcentration} is satisfied for a large class of distributions, including:

(a) $X_i$'s are Gaussian vectors (which is clearly a special case of \eqref{e.iid});

(b) $X_i$'s are sampled from the unit sphere.

(c) $X_i$'s are sampled from a distribution with log-concave density $e^{-U(x)}$ such that  $\text{Hess }U (x)- c(p)I_d$ is positive definite. (In this case $b=2$, see e.g.  \cite{ledoux}.)

Other examples can be found in \cite{karoui, karoui07a} and \cite{ledoux}.

In the special cases (a,b) above,  El Karoui's results were improved recently by Cheng--Singer \cite{cheng-singer}, where the authors showed that similar results hold for the variant \eqref{e.Adefnodiag} under the weaker assumption that $f$ is $C^1$ near $0$.

An initial examination of Karoui's results reveals that one only needs differentiability of $f$ at $0$ to formulate the above linear component $B$ of $A$. On closer looks, perhaps continuity of $f$ at $1$ is also required, since the diagonal entries of the covariance matrix of $X_i$ are converging to $1$ in the large $n$ large $p$ limit; except for the zero-diagonal model \eqref{e.Adefnodiag}. 

In the first result of this paper, it will be shown that under these minimal regularity assumptions on $f$ the nonlinear-to-linear results of \cite{karoui} can still be proved for a large class of distributions. Similar settings for distribution of $X_i$ will be considered: 
\begin{itemize}
\item The idd setting \eqref{e.iid} with $K>4$ moment bounds.
\item The high concentration setting \eqref{e.karouiconcentration} with  $c(p)\ge O(p^{b/2})$.
\end{itemize}
While our assumption $c(p)\ge O(p^{b/2})$ is stronger than El Karoui's assumption, it is in fact satisfied by a fairly large class of interesting distributions (see \cite{ledoux} and also \cite{karoui} for many examples); also a recent work of Gu\'edon--Milman \cite{guedon-milman} (cf. \cite{paouris}) indicates that such concentration assumption may be true for the isotrophic log-concave setting (see also the discussion following Conjecture~\ref{c.logconcave-inner} for details). On the other hand, Theorem~\ref{t.maincor1} requires only minimal regularity assumptions on $f$; one might view this as a trade-off between regularity of $f$ and concentration assumptions on the population distribution.

\begin{theorem}\label{t.maincor1} Assume the iid setting \eqref{e.iid} with $K>4$, or the high concentration setting \eqref{e.karouiconcentration} with $c(p)\ge O(p^{b/2})$.

(i) Let $f$ be differentiable at $0$ and continuous at $1$. Let $A$ be defined by \eqref{e.Adefdiag} with $F(X,Y,p)=f(X^TY)$. Then $A$ has the same limiting spectral distribution as
$$B = [f(1)-f(0)-f'(0)]I_n + f'(0) (X_i^TX_j)_{i,j} \ \ .$$

(ii) Let $f$ be differentiable at $0$. Let $A$ be defined by \eqref{e.Adefnodiag} with $F(X,Y,p)=f(X^TY)$. Then $A$ has the same limiting spectral distribution as $$B =  [-f(0)-f'(0)]I_n + f'(0) (X_i^TX_j)_{i,j} \ \ .$$
\end{theorem}

Theorem~\ref{t.maincor1} will be shown in Section~\ref{s.pindependent}.

\subsubsection{The distance kernel}

Let $F(X,Y,p)=f(\|X-Y\|^2)$. This model has recently attracted the attention of some authors 
(see e.g. \cite{karoui,  jiang, bordenave}), motivated by connections to machine learning theory and physics. In this model, it is clear that the two settings \eqref{e.Adefdiag} and \eqref{e.Adefnodiag} are equivalent up to a shift (by $f(0)$) of the limiting spectral distribution.  Below it will be assumed that $A$ is defined using \eqref{e.Adefdiag}.

When $(X_i)$ are Bernoulli or sampled from the unit sphere, the distance kernel model follows from the inner product model, but it is not hard to find  interesting examples (such as Gaussian or log-concave) where a naive adaptation of this reduction does not apply. This however suggests that  $A$ should have the same limiting spectral distribution as
$$B = [f(0) - f(2)+2f'(2)]I_n -2f'(2)(X_i^T X_j)_{i,j}$$
when $f$ is sufficiently smooth and the distribution of $X_i$ is sufficiently concentrated. This was shown in El Karoui \cite{karoui}, where the author assumed essentially the same settings for the $X_i$'s as in the last section:
\begin{itemize}
\item In the iid setting \eqref{e.iid} with $K>5$ moment bounds, this was shown for $f$  being $C^3$ near $2$.
\item In the high concentration setting \eqref{e.karouiconcentration} with $c(p)>O(p^{\epsilon+b/4})$ this was shown for $f$ being $C^2$ near $2$ and $C^1$ near $0$. 
\end{itemize}
In the iid setting, it was furthermore shown in \cite{karoui} that a stronger convergence holds in spectral norm, which may lead to more information about the distribution of the largest eigenvalue of $A$. (The interested reader is refered to to \cite{karoui, karoui07b} and refefences there-in for related literature.) As remarked earlier, the limiting spectral distribution of $B$ may be computed explicitly using Mar\v{c}henko--Pastur theory.

It is clear that one  only requires differentiability of $f$ at $2$ to write down the above linear component $B$ of $A$.  Theorem~\ref{t.maincor2} below shows that for a large class of distributions of $X_i$, El Karoui's nonlinear-to-linear results for distance random matrices can be proved for $f$ assuming only this differentiability. 

\begin{theorem}\label{t.maincor2} Assume the iid setting \eqref{e.iid} with $K>4$, or the high concentration setting \eqref{e.karouiconcentration} with $c(p) \ge O(p^{b/2})$.

Let $A$ be defined using  \eqref{e.Adefdiag} with $F(X,Y,p)  = f(\|X-Y\|^2)$ where $f$ is differentiable at $2$. Then $A$ has the same limiting spectral distribution as 
$$B = [f(0) - f(2)+2f'(2)]I_n - 2 f'(2) (X_i^T X_j)_{i,j} \ \ .$$
\end{theorem}

Besides regularity improvement,  in the iid setting Theorem~\ref{t.maincor2}  requires less  moment bounds than El Karoui \cite{karoui}. As discussed before in the paragraph leading to the statement of Theorem~\ref{t.maincor1}, our assumption on $c(p)$ is stronger than that in \cite{karoui}, but is satisfied by a large class of interesting distributions, see \cite{ledoux} and \cite{karoui} for many examples. A recent work of Gu\'edon--Milman \cite{guedon-milman} (cf. \cite{paouris}) indicates that such concentration inequality may hold in the isotrophic log-concave setting, see Conjecture~\ref{c.logconcave-distance} for details.

Theorem~\ref{t.maincor2} will be proved in Section~\ref{s.pindependent}. 
 
\subsubsection{Conjectures}

We conjecture  that similar results hold in the isotrophic log-concave case.  

\begin{conjecture}\label{c.logconcave-inner} The results of Theorem~\ref{t.maincor1} hold if  $X_i$'s are iid  random vectors from a normalized isotropic log-concave distribution. 
\end{conjecture}

\begin{conjecture}\label{c.logconcave-distance} The results of Theorem~\ref{t.maincor2} hold if  $X_i$'s are iid  random vectors from a normalized isotropic log-concave distribution. 
\end{conjecture}

Recall the following result of Gu\'edon-Milman \cite{guedon-milman} (cf. Paouris \cite{paouris}):
\begin{proposition} \label{p.guedo-milman} \emph{\cite{guedon-milman}}
Let $X$ be a (normalized) random vector sampled from an isotropic 
 measure in $\R^p$ with log-concave density, then for any $t \ge 0$
$$\P ( \big|\| X \|  -  1 \big| \ge t   ) \le C \exp(- c  p^{1/2} \min(t,t^3)) $$ for some absolute constants $c, C \in (0,\infty)$.
\end{proposition} 

It follows from the proofs of Theorem~\ref{t.maincor1} and Theorem~\ref{t.maincor2} that both conjectures hold if the estimate of Proposition~\ref{p.guedo-milman} could be improved to $O(\exp(-cp^{1/2}t))$. Mainly, in these proofs the concentration bound \eqref{e.karouiconcentration} is only needed for the Lipschitz function $\|X\|$, here $X$ is either $X_i$ or $X_i+X_j$ (which also have an isotrophic log-concave distribution).
This leads to the question of whether one can remove the term $t^3$ in the Gu\'edon-Milman concentration result.

See also Bordenave \cite{bordenave} for a recent work that was motivated by Conjecture~\ref{c.logconcave-distance}.

\subsubsection{Possible extensions of Theorem~\ref{t.maincor1} and Theorem~\ref{t.maincor2}}

Theorem~\ref{t.maincor1} and Theorem~\ref{t.maincor2} could be considered as extensions of Theorem 2.3 and Theorem 2.4 of El Karoui in \cite{karoui} when the envelope/kernel functions are rougher.  On the other hand, in \cite{karoui} El Karoui considered a more general setting when the covariance matrix $\Sigma$ of $X_i$ is less restrictive (in the current paper we assume $\Sigma=Id$ the identity operator in $\mathbb R^p$, which is the simplest but also most natural setting). More precisely, in \cite{karoui} $\Sigma$ is allowed to depend on $p$, but is still positive definite and converges in some fashion to identity in the limit $n,p\to\infty$. We anticipate that the proof of Theorem~\ref{t.maincor1} and Theorem~\ref{t.maincor2} in this paper has a natural extension that could lead to an extension of these Theorems to settings similar to those considered in Theorem 2.3 and 2.4 of \cite{karoui}, however these extensions are not explored in the current paper.

\subsection{The $p$-dependent setting}
In this section, a less classical setting recently investigated in \cite{cheng-singer} will be discussed. Here, the envelope function $f$ is allowed to depend on $p$ and may have {\it  very little regularity} 
 in $x$. In other words, $f$ may be varying with $p$. Examples of such situations and their motivations are presented in \cite{cheng-singer}. In this paper, only the inner-product kernel $g(X,Y)=X^T Y$ will be considered, and a similar investigation for the distance kernel is left for a further study.  Furthermore, following \cite{cheng-singer}, only the non-diagonal model \eqref{e.Adefnodiag} will be considered, and analogous results for the diagonal model \eqref{e.Adefdiag} may be obtained using a diagonal perturbation argument.  

In this section,  $X_1, X_2,\dots, X_n$ will be  iid random vectors in $\mathbb R^p$ whose coordinates are independent copies of a random variable $Z$ with mean $0$ and variance $1/p$, such that for all $K >0$ there is a constant $C_K$ depending on $K$ such that
\begin{equation}\label{e.momentassumption}
\E |Z|^K \le C_K \, p^{-K/2} \ \ . 
\end{equation}
While \eqref{e.momentassumption} requires all $K>0$, this assumption may be improved  if there are better bounds on the growth of a scaled version of $f$ as $p\to\infty$.  For details, see the remark after the statement of Theorem~\ref{t.main}.

Below, some standard facts about orthogonal polynomials will be recalled, for a standard reference see e.g. \cite{szego} or \cite{andrewsetal}. Given a nonnegative 
measure $\mu$ on $\R$ and $k=0,1,2,\dots$, the $k^{th}$ orthogonal polynomial $p_k(x)$ with respect to $\mu$ is a  polynomial of degree $k$ with positive leading coefficient, such that
$$\int_{\R} p_k(x) p_m(x) d\mu(x) = \begin{cases} 0, &\text{if $m\ne k$;}\\ 
1, &\text{if $m=k$.}\end{cases}$$
For any  function $h\in L^2(d\mu)$ that belongs to the span of $\{p_n, n\ge 0\}$, one has the formal series
$$\sum_{k=0}^\infty a_k \, p_k(x)\ \ , \qquad \text{where} \qquad a_k := \<h,p_k\>_{\mu} = \int_{\R} h(x) p_k(x) d\mu(x) \ \ ,$$
and if the series converges to $h$ in $L^2(\mu)$ then the Plancherel equality holds: $\|h\|_{L^2(d\mu)}^2 = \sum_{k\ge 0} |a_k|^2$.
In that case, since $p_0(x)\equiv 1$ for probability distributions, it follows that if $\mu$ is the distribution of some random variable $\xi$ then
$$\Var [h(\xi)] = \sum_{k\ge 1} |a_k|^2 \ \ .$$

For each $p$ let $\xi_p= \sqrt p X^T Y$ where $X$ and $Y$ are two iid copies of (any) vector $X_i$. Clearly $\xi_p$ has  mean $0$ and variance $1$. Let $p_{k,p}(x)$, $k\ge 0$, be the orthogonal polynomials with respect to the probability distribution $\mu_p$ of $\xi_p$. 

Below we state the conditions that will be assumed on the envelope function $f(x,p)$ for the next result, Theorem~\ref{t.main}. These conditions were first formulated in \cite{cheng-singer} in an equivalent form.  For any $f(x,p)$ let $k(x,p):= \sqrt p f(x/\sqrt p,p)$ and consider the expansion
$$k(x,p) \sim \sum_{k\ge 0} a_{k,p} p_{k,p}(x)  \ \ .$$	
In this paper, $f$ is said to be admissible with respect to the generating distribution of $X_i$ if the following three conditions hold.
\begin{itemize}
\item[(i)] \emph{(Uniform convergence)} The orthogonal polynomial series of $k(x,p)$ converges to $k(x,p)$ in $L^2(\mu_p)$ uniformly over $p$ large. In other words,  for any $\epsilon>0$ there exists $L=L(\epsilon)$ such that the following holds for $p$ large
\begin{equation}\label{e.kerneluniform} \|k-\sum_{i\le L} a_{i,p} p_{i,p}\|^2_{L^2(\mu_p)} = \sum_{i>L} |a_{i,p}|^2 \le \epsilon \ \ .
\end{equation}
\item[(ii)] \emph{(Normalization)} There exists $\nu \in [0,\infty)$ such that
\begin{eqnarray}
\label{e.kernelvariance}
\lim_{p\to\infty} \sum_{i=1}^\infty |a_{i,p}|^2 = \nu \ \ .
\end{eqnarray}
\item[(iii)] \emph{(Scaling)} There exists $a \in [0,\infty)$ such that
\begin{equation}\label{e.kernelscaling}
\lim_{p\to\infty} a_{1,p} = a\  \ .
\end{equation}
\end{itemize} 
It is clear that \eqref{e.kernelvariance} and \eqref{e.kernelscaling} together imply the condition $a^2 \le \nu$, which shall be assumed throughout.  It is worth pointing out that the set of orthogonal polynomials with respect to a probability measure $\mu$ does  not always form a complete basis in $L^2(\mu)$, this however holds for a fairly large class of probability measures, including those with sub-exponential tails (i.e. $\P(|\xi|>x) = O(e^{-c|x|})$ for some $c>0$, see e.g. \cite[Theorem 6.5.2]{andrewsetal}). In particular this completeness holds if the measure is compactly supported. It follows that the convergence of the orthogonal expansion in $L^2$ holds automatically if $X_i$'s are Gaussian or bounded (note that this gives convergence of the expansion for each $p$, and condition \eqref{e.kerneluniform} is about the uniformity of the convergence). In the general case when completeness of the orthogonal polynomials is not guaranteed, the condition \eqref{e.kerneluniform} has to be checked carefully (for both the convergence of the expansion for each $p$, and the uniformity of the convergence over $p$ large).

For convenience of the reader, the definition of the Stieltjes transform $m(z)$ of a measure $\mu$ is recalled below:
$$m(z) = \int_{\mathbb R} \frac{d\mu(x)}{x-z} \ \ , \ \ Im(z) > 0 \ \ .$$

\begin{theorem}\label{t.main} Assume that $f$ is admissible with respect to $X_i$'s which satisfy \eqref{e.momentassumption}. Let $A$ be generated using \eqref{e.Adefnodiag} using $F(X,Y,p) = f(X^T Y,p)$. Then, the empirical distribution of $A$ converges weakly to a probability distribution whose  Stieltjes  transform $m(z)$ satisfies:
\begin{equation}\label{e.functionaleqn}
-\frac 1 {m(z)} = z+ a\Big(1-\frac 1{1+\frac {a}{\gamma}m(z)}\Big) + \frac {\nu-a^2}{\gamma} m(z) \ \ .
\end{equation}
\end{theorem}

\begin{remark}
As pointed out in \cite{cheng-singer}, this limiting spectral distribution is no longer MP when $\nu\ne a^2$. Unique solvability of  \eqref{e.functionaleqn} was proved in \cite{cheng-singer} using elementary arguments. 
\end{remark}

The special case of Theorem~\ref{t.main} for Gaussian random vectors  was proved by Cheng--Singer in \cite{cheng-singer}.  Theorem \ref{t.main} 
 (positively) answered the question of Cheng and Singer in \cite{cheng-singer} about the validity of their result for more general distribution (such as Bernoulli).  
 
While it is assumed in Theorem~\ref{t.main} that \eqref{e.momentassumption} holds for all $K>0$,  the method of proof can be easily refined to lessen this assumption when more information about $f$ is given.  More precisely, if there is an upper bound $L$ on the cut-off degree $L(\epsilon)$ in \eqref{e.kerneluniform} (independent of $\epsilon>0$) then we only need to have \eqref{e.momentassumption} up to $O_L(1)$: the degree $L$ will enter the proof in Lemma~\ref{l.entrydecomposition} and Lemma~\ref{l.momentbound} and will eventually dictate the number of required moment bounds on entries of $X_i$'s. For instance, when $f$ is independent of $p$  and  has a non-vanishing derivative at $x=0$ it can be verified that $\nu=a^2=f'(0)^2$, and one can see from the proof of Theorem~\ref{t.maincor1} that one could take $L=1$ in \eqref{e.kerneluniform}. Eventually, with some refinements (tailored specifically for Theorem~\ref{t.maincor1}), this leads to the requirement $K>4$ in Theorem~\ref{t.maincor1}.

{\it Acknowledgement.} We would like to thank  X. Cheng and  A. Singer for bringing this interesting subject to our attention and many useful conversations. We would like to thank the referees for corrections and  suggestions which have lead to improvement of the quality of the paper.

\section{The general ideas}

Let $m_A(z)$ denote the Stieltjes transform of the empirical spectral distribution $\rho_A$ of $A$; in the following $m_A$ will be refered to as the Stieltjes transform of $A$. Explicitly,
$$m_A(z) = \frac 1 n \Tr[(A-z)^{-1}] \ \ ,  \ \ \text{Im}(z)>0.$$
By standard reductions (see e.g. \cite{bai-silverstein}), it suffices to show that
 $m_A(z) $ converges  to the Stieltjes transform of the desired limiting spectral distribution (which is always a probability distribution in the current paper) for $Im(z) >0$. For instance, in the setting of Theorem~\ref{t.main} it will be shown that $m_A$ converges to the solution of \eqref{e.functionaleqn}. 

The main idea for showing  the desired convergence of $m_A$ is to compare  $A$ with a suitably chosen random matrix whose Stieltjes transform
 already has the desired convergence.   In fact, due to a result from   \cite{cheng-singer}
 which asserts that 
\begin{equation}\label{e.concentration}
\lim_{n\to\infty} |m_A(z) -\E m_A(z)| = 0 \ \ \text{a.s.},
\end{equation}
it suffices to compare expected values of Stieltjes transforms in question. 
   To keep the paper self-contained, a short  proof of \eqref{e.concentration} will be included in section~\ref{s.concentration}  (see Lemma~\ref{l.concentration}).
   
   The proof for the results in the $p$-independent case will be presented in the next section. Following El Karoui \cite{karoui},  $A$ will be compared with a linear approximation of $A$, obtained by replacing the envelope function $f$ with the linear part of its Taylor expansion at a suitable point.  The main idea which allows us to improve the regularity assumptions on $f$ in El Karoui and Cheng--Singer's results is a simple \emph{transference principle}, see Lemma~\ref{l.transference} and also its companion Lemma~\ref{l.removediagonal} in Section~\ref{s.pindependent} for details.

The proof of the results in the $p$-dependent setting will use a series of comparisons.
In order to carry out the  analysis of the main comparison, the Lindeberg 
swapping method will be used, following ideas from \cite{chatterjee, taovu12acta}. 
This method has recently proved useful in various studies of random matrices, especially 
for the local statistics (see \cite{taovusurvey} for a survey). One of the main difficulty in implementing the Lindeberg method in the current setting is the lack of regularity of the envelope function. To overcome this difficulty, the uniform convergence condition \eqref{e.kerneluniform}  will be used, as this condition allows for approximation of $f$ with polynomials (which are very smooth). The proof of Theorem \ref{t.main} will be presented in Section~\ref{proof1} and Section~\ref{s.main}.

In the rest of the paper,  without loss of generality it will be assumed that $\text{Im}(z)>0$. All implicit constants in the paper may depend on $z$.
All asymptotics notations  are used under the assumption that $p,n \rightarrow \infty$.

\section{The $p$-independent setting}\label{s.pindependent}

In this section, we prove Theorem~\ref{t.maincor1} and Theorem~\ref{t.maincor2}.

Let $A$ be defined by \eqref{e.Adefnodiag} using $f(g(X,Y))$ and let $\mathcal A$ be defined by \eqref{e.Adefnodiag} using $F(X,Y,p)=g(X,Y)$. The following {\it  transference principle} asserts that one can deduce the limiting spectral density for  $A$ from $\mathcal A$ as long as:

\indent (i) $f$ is differentiable at  the mean value of $g(X_i,X_j)$; and

\indent (ii) the entries of $\mathcal A$ (hence the kernel $g$) satisfies a fairly general concentration condition (relative to $X_i$'s).

\begin{lemma}[Transference principle]  \label{l.transference}
Assume that $\rho_{\mathcal A}$ converges weakly to a probablilty distribution.  
Let $a=\E g(X_i,X_j)$ and let $f$ be differentiable at $x=a$. Assume
\begin{equation}\label{e.entrynormalization}
\Var [g (X_i,X_j)] = O(1/p)  \ \ , 
\end{equation}
and assume  that for any fixed $\delta>0$ it holds that
\begin{equation}\label{e.entryconcentration}
\P ( \max_{i \neq j}  | g(X_i, X_j) -a | > \delta) =o(1)   \ \ .
\end{equation}
Then $A$ has the same limiting spectral distribution as $$B = (af'(a)-f(a)) I_n + f'(a)\mathcal A \ \ .$$

\end{lemma} 

Remarks:  While different pairs $(f,g)$ may generate the same $F=f\circ g$, the two constraints \eqref{e.entrynormalization} and \eqref{e.entryconcentration}  impose a strong normalization on $g$. Also, in Lemma~\ref{l.transference} the spectral distribution of $\mathcal A$ is not required to be  Mar\v{c}henko--Pastur.

The following simple result will also be used, which says that under an assumption on concentration of $g(X_i,X_i)$,  the models \eqref{e.Adefdiag} and \eqref{e.Adefnodiag}  are equivalent. 

\begin{lemma}   \label{l.removediagonal} Let  $A_1$ and $A_2$ be defined by \eqref{e.Adefdiag} and \eqref{e.Adefnodiag} respectively using $F(X,Y,p) = f(g(X,Y))$. Assume that $f$ is continuous at $b:=\E g(X_i,X_i)$. Assume that for any $\delta>0$ it holds that
\begin{equation}\label{e.diagonalconcentration}
P (\max_{1\le i \le n} \big|g(X_i,X_i)-b\big|>\delta) = o(1)  \ \ .
\end{equation}
Then in the large $n$ large $p$ limit it holds that
$$|m_{A_1}(z+f(b))-m_{A_2}(z)| = o(1) \ \ \text{a.s.}$$
\end{lemma}
\begin{remark} If $g(X_i,X_i)$ is a constant then the above result is trivial, in which case  continuity of $f$ at $b$ is not needed. 
\end{remark}

Using Lemma~\ref{l.removediagonal}, the main argument is reduced to the non-diagonal model \eqref{e.Adefnodiag}, where the transference principle could be used. Lemma~\ref{l.removediagonal} could be viewed as a companion of Lemma~\ref{l.transference}.

Below, Theorem~\ref{t.maincor1} and Theorem~\ref{t.maincor2} are deduced from the above two Lemmas. Proofs of Lemma~\ref{l.transference}  and Lemma~\ref{l.removediagonal} are presented in Section~\ref{s.transference}.

\subsection{Proof of Theorem~\ref{t.maincor1} }

\noindent \underline{The iid case}: Assume that the entries of $X_i$ are iid with $K>4$ moment bounds. 

\emph{Step 1:} We first reduce the Theorem to the model \eqref{e.Adefnodiag} of $A$. By Lemma~\ref{l.removediagonal}, it suffices to show that: for any fixed $\delta>0$ (i.e. independent of $n,p$) and any $i$ we have
\begin{equation}\label{e.verifylengthconcen}
P( \big|\|X_i\|^2-1\big|>\delta) = o(1/p) \ \ ,
\end{equation}
We note that a more quantative estimate was proved in \cite{karoui} using more careful arguments, on the other hand for \eqref{e.verifylengthconcen} the following simplified argument suffices. Fix $\delta>0$ and let $X=(x_1,\dots, x_p)$ be an independent copy of $X_i$'s. Let $M:=p^\beta$ for $0> \beta> \frac 2 K - \frac 1 2$. Let $E=\{\max_{i} |x_{i}| > M\}$, clearly
$$P(E) \le   C \sum_{j} p^{-\beta K} \E |x_{j}|^K = o(1/p) \ \ .$$
Let $\widetilde X = (x_j 1_{|x_j|\le M})_{j=1}^p$. On $E^c$ clearly $\widetilde X=X$. Thus, it  suffices to show that
\begin{equation}\label{e.rowconcen}
P(\big|\|\widetilde X\|^2 - 1\big|>\delta) = o(1/p^2)  \ \ .
\end{equation}
Let $\vec 1 = (1,\dots , 1) \in \R^p$. Let $\mu$ and $\sigma^2$ be the   mean and variance of  $x_i1_{|x_i|\le M}$.  It is not hard to see that
$$\mu=  o(\frac 1 {\sqrt p}) \qquad \text{and}\qquad \sigma^2 = \frac 1 p + o(\frac1 p) \ \ .$$
For $p$ sufficiently large, it follows that
\begin{equation}\label{e.rowrenormalize}
P(\big|\|\widetilde X\|^2 - 1\big|>\delta)  \le P(\big|\|\widetilde X - \mu \vec 1\|^2 - p\sigma^2\big|>\delta/2) \ \ .
\end{equation}
Write $\|\widetilde X - \mu \vec 1 \|^2 - p\sigma^2 = \sum_{j=1}^p \big[ (\widetilde x_j - \mu)^2 -\sigma^2\big]$  sum of iid random variables, each has mean $0$  and is bounded above by $O(M^2)=O(p^{2\beta})$ and has the following variance bound:
$$N:=\Var[(\widetilde x_j - \mu)^2 -\sigma^2] = \E[ (\widetilde x_j-\mu)^4 ]- \sigma^4 = O(p^{-2}) \ \ .$$
By Chernoff's inequality (see e.g. \cite{taobook}), for $C_1, C_2$ absolute positive constants it holds that
$$P(\big|\|\widetilde X - \mu \vec 1\|^2 - p\sigma^2\big|>\delta/2)  \le $$
\begin{equation}\label{e.truncatedrowconcen}
\le C_1 \max \Big(\exp(-\frac{C_2 \delta^2}{p N}), \exp(-\frac{C_2 \delta}{M^2})\Big) = o(p^{-1}) \ \ .
\end{equation}
Collecting inequalities \eqref{e.rowrenormalize}, and \eqref{e.truncatedrowconcen}, the desired estimate \eqref{e.rowconcen} follows.

\emph{Step 2:} Thanks to Step 1, it remains to show the theorem for $A$ given by \eqref{e.Adefnodiag}. Let $\delta>0$ be fixed. Using Lemma~\ref{l.transference}, it suffices to show that
\begin{equation}\label{e.innerproduct}
\E\big[|X_i^T X_j|^K\big] = O(p^{-K/2}) \ \ , \ \   \text{for some $K>4$.}
\end{equation}
Write $X_i = (x_{i1},\dots, x_{ip})$ and $X_j=(x_{j1},\dots, x_{jp})$.  
Then $x_{ik}$ and $x_{jm}$ are independent with mean $0$ and variance $1/p$ for any $1 \le k, m \le p$. 
By the inverse Khintchine inequality (i.e. the Marcinkiewicz--Zygmund inequality) it holds that
\begin{eqnarray*}
&&\E |X_i^T X_j|^{K} \le C_K \E( \sum_{m=1}^p |x_{im}  x_{jm} |^2)^{K/2}\\
&\le& C_K p^{K/2-1} \sum_{m=1}^p \E |x_{im} |^K \E |x_{jm}|^K  \qquad \text{(H\"older, then independence)}\\
&=& O(p^{-K/2}) \qquad \text{(using given moment bounds).}
\end{eqnarray*}

\underline{The high concentration case}:
As in the iid case, it suffices to show \eqref{e.verifylengthconcen} and \eqref{e.innerproduct}.  As can be seen below, in the proof it is enough to assume \eqref{e.karouiconcentration} for the $1$-Lipschitz functions of the form $f(X)= \|X+c\|$, $c\in \R^p$ constant vectors. See also the discussion after the statements of Conjecture~\ref{c.logconcave-inner} and Conjecture~\ref{c.logconcave-distance}. 

\emph{Proof of \eqref{e.verifylengthconcen}:} Let $F(Y)=\|Y\|$ the Euclidean length of $Y\in \R^p$, clearly $F$ is $1$-Lipschitz. Fix $i$ and $\delta>0$, without loss of generality assume $\delta<1/2$. 

We first show that uniformly over $r>0$ it holds that
\begin{equation}\label{e.concen1}
P(\big|\|X_i\|-1\big| > r)  = O(e^{-c(p)r^b/C_b}) \ \ ,
\end{equation}
here and below $C_b$ will denote absolute constants that could depend on $b$. Let $a=\E \|X_i\|$. It sufficies to show that $|a-1| = O(c(p)^{-1/b})$.  By \eqref{e.karouiconcentration}, it holds that
$$|a-m_F| \le   \E \big|\|X_i\|-m_F\big| $$
$$= O(\int_0^\infty e^{-c(p)r^b}dr ) = O(c(p)^{-1/b}) \ \ .$$
Let $C$ be the implicit constant in the last estimate. Then for any $r>2C c(p)^{-1/b}$ it holds that
$$P(\big|\|X_i\|-a\big| > r)  = O(e^{-c(p)(r/2)^b}) \ \ .$$
In this estimate, it is clear that if $r = O(c(p)^{-1/b})$ then $e^{-c(p)(r/2)^b} \sim 1$ while the left hand side is at most $1$. Thus, the above estimate holds uniformly over $r>0$.  Now, \eqref{e.concen1} follows from
$$0 \le \E \|X_i\|^2 -a^2 = \E \big|\|X_i\|-a\big| ^2 = O(c(p)^{-2/b}) \ \ .$$

We obtain
\begin{equation}\label{e.Xisquare}
P(\big|\|X_i\|^2-1\big| > r) = \begin{cases} O(e^{-c(p)r^b/C_b}), & \text{if $r\ge O(1)$;}\\
O(e^{-c(p)r^{b/2}/C_b}), & \text{if $r=O(1)$.}
\end{cases}
\end{equation}
In particular,  \eqref{e.verifylengthconcen} follows.

\emph{Proof of  \eqref{e.innerproduct}:}  We first show that
\begin{equation}\label{e.concentrationXiXj}
P(\big|\|X_i+X_j\|- \E\|X_i+X_j\|\big|>2r) = O(e^{-c(p) r^b}) \ \ , \ \  i\ne j \ \ ,
\end{equation}
uniformly over $r>0$. For any $X\in \R^p$ let $G(X) = \E_{Y}\|X+Y\|$ the expectation over $Y$ independently sampled from the distribution of $X_i$'s. It is clear that
$$\E_{X_i} G(X_i) = \E_{X_i,X_j}\|X_i+X_j\|  \ \ . $$
Using  independence of $X_i,X_j$, it follows that
$$\text{LHS of \eqref{e.concentrationXiXj}} $$
$$\le  \E_{X_j} [\E_{X_i} [1_{\big|\|X_i+X_j\|- G(X_i)\big| > r}]] +  \E_{X_i} [1_{\big|G(X_i) - \E_{X_i} G(X_i)\big| > r}]$$
$$=  \E_{X_j} O(e^{-c(p)r^b}) + O(e^{-c(p)r^b})$$
therefore \eqref{e.concentrationXiXj} follows. 

We now show that
\begin{equation}\label{e.Xi+Xj}
P(\big|\|X_i+X_j\|-\sqrt 2\big|> r) = O(e^{-c(p)r^b/C_b}) \ \ ,
\end{equation}
uniformly over $r>0$. Let $\alpha = \E\|X_i+X_j\|$. It suffices to show that $\alpha=\sqrt 2 + O(c(p)^{-1/b})$.
For any $K\ge 1$ we have
$$\E \big|\|X_i+X_j\|- \alpha\big|^K =  \int_0^\infty K r^{K-1} P (\big|\|X_i+X_j\|- \alpha\big|>r) dr$$
$$= O(c(p)^{-K/b})\ \ .$$
Letting $K=2$ in the above estimate, it follows that $\alpha=\sqrt 2 + O(c(p)^{-2/b})$.

It follows from \eqref{e.Xi+Xj} that
$$P(\big|\|X_i+X_j\|^2- 2\big|>r)   = \begin{cases} O(e^{-c(p)r^{b/2}/C_b}) , &r \ge O(1) \ \ ;\\
O(e^{-c(p)r^{b}/C_b}) ,  & r = O(1) \ \ .
\end{cases}$$
Combining this with \eqref{e.Xisquare}, it follows that
\begin{equation}\label{e.XiXj}
P(|X_i^T X_j| > r)   = \begin{cases} O(e^{-c(p)r^{b/2}/C_b}) , &r \ge O(1) \ \ ;\\
O(e^{-c(p)r^{b}/C_b}) ,  & r = O(1) \ \ .
\end{cases}
\end{equation}
Consequently, for any $K\ge 1$ it holds that
$$\E [|X_i^T X_j|^K] = \int_0^{\infty} Kr^{K-1}P(|X_i^T X_j| > r) dr $$
$$= \int_0^1 + \int_1^\infty = O(c(p)^{-K/b}) + O(c(p)^{-2K/b}) = O(p^{-K/2}) \ \ .$$

\subsection{Proof of Theorem \ref{t.maincor2} }
For the distance model, the diagonal entries are $f(0)$, therefore removing/adding these entries does not require any regularity of $f$. The transference principle Lemma~\ref{l.transference} will be used, and it remains to show that $g(X,Y)=\|X-Y\|^2$ satisfies the two kernel conditions of Lemma~\ref{l.transference}. 

\underline{The iid case:} We first verify \eqref{e.entrynormalization}. Let $X$ and $Y$ denote $X_i$ and $X_j$ for some $i\ne j$. Then
$$\Var [\|X- Y\|^2]  = -4 + \sum_{i=1}^p \E (x_i -y_i)^4 + \sum_{i \neq j} \E (x_i- y_i)^2 (x_j-y_j)^2 $$
Using $\E x_i^4$ and $\E y_i^4 =O( 1/p^2)$, it is clear that
$$\sum_{i=1}^p \E (x_i -y_i)^4 = O(1/p) \ \ .$$
Using also independence of $x_i,x_j,y_i,y_j$, we have
$$ \sum_{i \neq j} \E (x_i- y_i)^2 (x_j-y_j)^2 = 4\sum_{k} \E x_k^2 \sum_{k} \E y_k^2 + O(1/p) = 4 + O(1/p) \ \ , $$
and \eqref{e.entrynormalization} follows.

We now verify  \eqref{e.entryconcentration}. Fix any $\delta>0$. Using the previously obtained bounds \eqref{e.verifylengthconcen}  and \eqref{e.innerproduct}, it follows from the triangle inequality that
$$\P( \max_{i\ne j} \big|\|X_i-X_j\|^2 -2\big| \ge \delta) $$
$$\le 2 P( \max_{i} \big|\|X_i\|^2 - 1 \big| > \delta/4) + P(\max_{i\ne j} |X_i^T X_j| > \delta/2)$$
$$= o(1) \ \ ,$$
as desired.

\underline{The high concentration case:} Note that \eqref{e.entryconcentration} follows from \eqref{e.Xisquare} and \eqref{e.XiXj}. In fact, for any $r>0$ it holds that
$$P(\big|\|X_i-X_j\|^2-2\big|>3r)\le 2 P(\big|\|X_i\|^2-1\big|>r) + P(|X_i^T X_j| > r)$$
$$=  \begin{cases} O(e^{-c(p)r^{b/2}/C}) , &r > 2 \ \ ;\\
O(e^{-c(p)r^{b}/C}) ,  & r = O(1) \ \ .
\end{cases}$$
Via the same argument as before, we also obtain
$$\E \big|\|X_i-X_j\|^2-2\big|^K = O(p^{-K/2}) \ \ ,$$
for any $K>2$, and taking $K=2$ gives us the first kernel condition.

\subsection{Proof of Lemma~\ref{l.removediagonal} and Lemma~\ref{l.transference}} \label{s.transference}
\proof[Proof of Lemma~\ref{l.removediagonal}]
Using \eqref{e.concentration}, it suffices to show that $\E |m_{A_1}(z+f(b)) - m_{A_2}(z)| = O(\epsilon)$ for any $\epsilon>0$, which is fixed in the rest of the proof.  

Let $\mathcal I_n$ denote the $n\times n$ identity  matrix. By a standard argument,
$$|m_{A_1}(z+f(b))-m_{A_2}(z)| \le $$
$$ \le \|(A_1-(z+f(b))\mathcal I_n)^{-1} - (A_2-z \mathcal I_n)^{-1}\|$$
here $\|.\|$ denotes the spectral norm of a matrix,
$$ \le C\|A_1-f(b){\mathcal I}_n-A_2\|$$
$$\le C\sup_{1\le i \le n} |f(g(X_i,X_i)) - f(b)| \ \ . $$
Since $f$ is continuous at $b$, there exists $\delta=\delta(f,\epsilon)>0$  such that 
$$|f(x)-f(b)|<\epsilon \ \ \text{if} \ \  |x-b|\le \delta \ \ . $$
Since $m_{A_1}(z), m_{A_2}(z)=O(1)$, it follows that
$$\E |m_{A_1}(z+f(b))-m_{A_2}(z)| =$$
$$= O(\epsilon) + O\Big(P(\sup_{1\le i \le n} |f(g(X_i,X_i)) - f(b)| > \delta) \Big) \ \ .$$
Therefore $\E |m_{A_1}(z+f(b))-m_{A_2}(z)|  = O(\epsilon)$, thanks to  \eqref{e.diagonalconcentration}.
\endproof

\proof[Proof of Lemma~\ref{l.transference}]

Without loss of generality we may assume that $a=0$. Recall that $A$ is defined using \eqref{e.Adefnodiag}.

Let $h(x)=f(0)+f'(0)x$, and let $B$ be obtained from $A$ by replacing $f$ with $h$ (while keeping the same kernel $g$). More specifically,
$$B = f(0) \mathcal M_1 + f'(0) \mathcal A - f(0) {\mathcal I}_n$$
here ${\mathcal I}_n$ is the $n\times n$ identity matrix and $\mathcal M_1$ is the  $n\times n$ matrix whose entries are all $1$'s (in particular $\mathcal M_1$ has rank $1$ and thus does not contribute to the limiting spectral distribution of $B$, see e.g. \cite{bai-silverstein} or \cite{bai}). Thus, it suffices to show that $m_A(z) - m_B(z)\to 0$ for any $z$ in the upper half plane.  Using \eqref{e.concentration}, this follows from
\begin{equation}\label{e.coupling1}
 |\E m_A(z) - \E m_B(z)|  \to  0  
\end{equation}
which will be shown in the rest of the proof.

Fix $\epsilon>0$, it suffices to show that $\E |m_A(z)-m_B(z)| =O(\epsilon)$ for $n$, $p$ sufficiently large. Let $\lambda_1(A) \le \dots \le \lambda_n(A)$ be the eigenvalues of $A$ and $\lambda_1(B)\le \dots \le \lambda_n(B)$ be the eigenvalues of $B$. Then for any fixed $z\not\in \mathbb R$ it holds that
\begin{eqnarray}
\nonumber &&|m_A(z) - m_B(z)|^2\\
\nonumber  &\le&  \frac 1 n \sum_{i=1}^n |\frac 1 {\lambda_i(A)-z} - \frac 1 {\lambda_i(B)-z}|^2 \qquad \text{(Cauchy-Schwarz)}\\
\label{e.pre-spectralnorm} &\le& C\frac 1 n \sum_{i=1}^n |\lambda_i(A)-\lambda_i(B)|^2 \qquad \text{(using $Im(z)>0$)}\ \ ,
\end{eqnarray} where $C >0$ is a constant which may depend on $z$. 
It then follows from the Hoffman--Wielandt inequality (see e.g. \cite{taobook})  that
\begin{equation} \label{e.pointwise}
|m_A(z) - m_B(z)|^2  \le C \frac 1 n \sum_{i\ne j}|A_{ij} - B_{ij}|^2 .
\end{equation}

By definition, there is $\delta>0$ depending on $f$ only such that 
\begin{equation}\label{e.differentiable}
|f(x)-h(x)|\le \epsilon |x| \qquad \text{for $|x|<\delta$.}
\end{equation} 
Let $F$ be the event that there is a pair $i \neq j$ such that $|g(X_i, X_j)| \ge \delta$. It follows from the second assumption \eqref{e.entryconcentration} on $g$ that  $\P(F) \le \epsilon^2$ for $n$ (and $p$) sufficiently large. We now estimate
$$\E [|m_A(z)-m_ B(z)|^2] \le \E [1_F |m_A(z)-m_ B(z)|^2] +  \E [1_{F^c} |m_A(z)-m_ B(z)|^2] \ \ .$$
Since $m_A(z) = O(1)$ and $m_B(z)=O(1)$, it follows that
$$\E [1_F |m_A(z)-m_ B(z)|^2] = O(\P(F))  = O(\epsilon^2) $$
for large $n$ large $p$. On the other hand, from \eqref{e.pointwise} and \eqref{e.differentiable}, it follows that
$$\E [1_{F^c} |m_A(z)-m_B(z)|^2 ]= O( \frac 1 n \sum_{i\ne j} \E (\epsilon^2 |g(X_i, X_j)|^2))  \ \ .$$
With the first assumption on $g$, it follows that
$$\E [1_{F^c} |m_A(z)-m_B(z)|^2 ]= O( \epsilon^2) \ \ .$$
Consequently, in the large $n$ large $p$ limit it holds that
$$\E |m_A(z)-m_B(z)| \le (\E |m_A(z)-m_B(z)|^2)^{1/2} = O(\epsilon) \ \ .$$
This completes the proof of \eqref{e.coupling1}.

\endproof

\subsection{Concentration of the Stieltjes transform}\label{s.concentration}

In this section, we include a proof of   \eqref{e.concentration}. 
\begin{lemma}\label{l.concentration} \emph{\cite{cheng-singer}} Let $M$ be an $p\times n$ matrix with independent entries.
Let $A$ be defined by $A_{ij}=F(M_i, M_j,p)$ any real-valued function that is symmetric in the first two variables, here $M_i$ denotes the $i^{th}$ column of $M$. Let $\widetilde A$ have the same non-diagonal entries as $A$, and zero diagonal entries. Fix $z$ with $Im(z)>0$. Then
\begin{equation}\label{e.mAconcentration1}
\lim_{n\to \infty} |m_A(z) - \E m_A(z)| = 0 \ \ \text{a.s.}
\end{equation}
\begin{equation}\label{e.mAconcentration2}
\lim_{n\to \infty} |m_{\widetilde A}(z) - \E m_{\widetilde A}(z)| = 0 \ \ \text{a.s.}
\end{equation}
\end{lemma}

\proof The proof  largely follows the argument in \cite{cheng-singer}, which is a variant of standard arguments (see e.g. \cite{bai-silverstein}). Only \eqref{e.mAconcentration1} will be proved below, the proof for \eqref{e.mAconcentration2} is entirely similar. Using Borel-Cantelli's lemma, it suffices to show
\begin{equation}\label{e.borelcantelli}
\E |m_A(z) - \E m_A(z)|^4 = O(n^{-2}) \ \ .
\end{equation}
Recall the following standard estimate (i.e. Khintchine's inequality) for the martingale square function:
\begin{equation}\label{e.burkholder}
\E |\sum_{j=1}^n \Delta_j |^p \le C_p \E (\sum_{j=1}^n |\Delta_j|^2)^{p/2} \ \ , \ \ 1<p<\infty \ \ ,
\end{equation}
where $\Delta_j = S_{j+1}-S_j$ the martingale difference sequence. Below, a martingale will be constructed such that $S_0=\E m_A(z)$ and $S_n=m_A(z)$, and then show that the corresponding right hand side of \eqref{e.burkholder} is bounded above by a suitable power of $n$.

For any $0\le k\le n$ let $\cala_k$ be the sigma algebra generated by the last $k$ columns of $M$. Then define $S_k:=  \E_k \Tr(A-z)^{-1}$ which is a martingale with respect to the filtration $\{\cala_0 \subset \dots \subset \cala_n\}$ . Since $m_A(z)$ is measurable with respect to $\cala_n$, the construction gives $S_n =  n m_A(z)$ while clearly $S_0 = n\E m_A(z)$. It follows from \eqref{e.burkholder}, with $p=4$, that
$$\E |m_A(z)-\E m_A(z)|^4  \le \frac C{n^2}\, \E \Big( \max_{1\le j\le n} |\Delta_j|^4 \Big)\ \ .$$
It therefore suffices to show that,  uniformly over $j$,
\begin{equation}\label{e.martingalediff}
\E_j [\Tr(A-z)^{-1}]-\E_{j-1} [\Tr (A-z)^{-1}] = O(1)
\end{equation}
where the implicit constant is allowed to depend on $z$. Fix $1\le j \le n$ and let $B$ denote the $(n-1)\times (n-1)$ submatrix of $A$, obtained by deleting the $j^{th}$ row and the $j^{th}$ column of $A$. By definition of $A$, the entries of $B$ are independent of $M_j$, therefore $\E_j B = \E_{j-1} B$ and consequently it suffices for \eqref{e.martingalediff} to show that
$$|\Tr(A-z)^{-1} - \Tr(B-z)^{-1}|   = O(1) \ \ .$$
Since $f$ is symmetric and real-valued, the eigenvalues of $A$ and $B$ are real valued. Furthermore, they interlace by the Cauchy interlacing theorem (see e.g. \cite{taobook}). The desired estimate now follows immediately:
$$|\Tr(A-z)^{-1} - \Tr(B-z)^{-1}| \le \int |\frac d {dt} (\frac 1{t-z})| dt =O(\frac 1 {Im(z)}) \ \ .$$
\endproof

\section{The $p$-dependent setting}

\subsection{Some estimates for orthogonal polynomials} \label{proof1} 

In this section, some basic estimates involving orthogonal polynomials are proved, and these estimates will be used in the proof of Theorem~\ref{t.main}. 

Let $h_k(x)$ denote the $k^{th}$ Hermite polynomial, i.e. the orthogonal polynomial with respect to the Gaussian measure $\mu(x) = \frac 1 {\sqrt {2\pi}} e^{-x^2/2}$.   

Let $p_{k,p}$ denote the $k^{th}$ orthogonal polynomial with respect to $\mu_p$ the probability distribution of $\xi_p=\sqrt p Z_1^T Z_2$ where $Z_1, Z_2$ are iid random $1\times p$ vectors whose coordinates are independently sampled from a random variable $Z$ satisfying 
$$\E [Z] = 0 \ \ , \ \ \Var [Z] = 1/p \ \ , $$
and for some $K>0$ sufficiently large
\begin{equation}\label{e.K}\ \  \E \big[|Z|^K\big] = O(p^{-K/2}) \ \ .
\end{equation}

\begin{lemma}\label{l.opcompare} 
Let  $k\ge 0$. If \eqref{e.K} holds for $K$ sufficiently large then for any $\delta>0$ it holds that
\begin{equation}\label{e.opcompare}
|p_{k,p}(x)- h_k(x)| \le C\delta(1+|x|^k)  \ \ , \ \ x\in \R \ \ ,
\end{equation}
here the implicit constant $C$ may depend on $k$ (but not on $p$).
\end{lemma}

\proof For $k\ge 0$ let $m_k = \E \xi_p^k$. Then for some normalization  constant $c_k$ it holds that (see e.g. \cite{szego}):
$$p_{k,p}(x) = c_k \det \begin{pmatrix} 
m_0 & m_1 & \dots & m_k\\
        &         & \dots & \\
m_{k-1} & m_{k} & \dots & m_{2k-1}\\
1      &    x   & \dots  & x^k
\end{pmatrix} \ \ .$$
The leading coefficient of $p_{k,p}$ is $c_k \det M_{k-1}$ where
$$M_j =\begin{pmatrix} 
m_0  & \dots & m_j\\
        & \dots & \\
m_j   & \dots & m_{2j}
\end{pmatrix} \ \ .$$
Since $\|p_{k,p}\|_{L^2(\mu_p)} = 1$, it follows that
$$c_k^2 = \frac{1}{(\det M_{k-1})(\det M_k)} \ \ .$$
and the sign of $c_k$ is completely determined from the sign of $\det M_{k-1}$. Therefore in order to show \eqref{e.opcompare} it suffices to show that if $\delta>0$ then
\begin{equation}\label{e.momentcompare}
\E \xi_p^k = \E N^k + O_k(\delta)  \  \  ,
\end{equation}
where $N$ is the normal Gaussian $\mathcal{N}(0,1)$. But this is a classical theorem of von Bahr \cite{vonbahr}.
\endproof

\begin{lemma}\label{l.opcoefficients} Let $a_{j,p}$ be the coefficients in the orthogonal expansion of a normalized kernel function $k$ satisfying $\E [k(\xi_p,p)^2] = O(1)$ uniformly over $p$. Then
$$|a_{j,p}| = O(1) \ \ .$$
\end{lemma}
\proof  Clearly for any $j$
$$|a_{j,p}| = |\int k(x,p) p_{j,p}(x) d\mu_p(x)|$$
$$\le \|k(.,p)\|_{L^2(\mu_p)} \|p_{j,p}\|_{L^2(\mu_p)} = O(1) \ \ .$$
\endproof
Note that if the kernel $k$ satisfies \eqref{e.kernelvariance} then  Lemma~\ref{l.opcoefficients} applies.

\subsection{Proof of Theorem~\ref{t.main}}\label{s.main}

Let $M_A$ be the $p\times n$ matrix whose columns are $X_1,\dots X_n$. Let $M_G$ denote the $p\times n$ matrix whose columns are iid random Gaussian vectors $G_1$, ..., $G_n$, which are normalized so that the entries of $M_G$ has mean $0$ and variance $1/p$. A Gaussian analogue $G$ of $A$ will be constructed as follows. The construction will ensure, thanks to \cite{cheng-singer}, that the spectral density of $G$ converges to the desired limiting spectral density in Theorem~\ref{t.main}.

Let $\xi_{G,p}=\sqrt p G_1^T G_2$ and let  $\mu_{G,p}$ denote its probability distribution. Let
$$K(x,p) = \sum_{i=0}^\infty a_{i,p} P_{i,p}(x)$$ 
where $P_{i,p}$ is the $i^{th}$ orthogonal polynomial with respect to $\mu_{G,p}$. Note that the above infinite sum converges uniformly in $L^2(\mu_{G,p})$ by given assumptions on $a_{i,p}$. 

Let $F(x,p)=(\sqrt p)^{-1} K(x\sqrt p, p)$ and let $G$ be the random matrix generated from $M_G$ using $F$. It was shown in \cite{cheng-singer} that the Stieltjes transform $m_G(z)$ of $G$ converges pointwise to the solution of \eqref{e.functionaleqn}. Therefore, using \eqref{e.concentration}, it suffices to show that in the large $n$ large $p$ limit it holds that
\begin{equation}\label{e.coupling2}
\E [m_{A}(z)  - m_G(z)] =O(\epsilon) \ \ ,
\end{equation}
where $\epsilon>0$ is fixed in the rest of this section.

It follows from the given assumption \eqref{e.kerneluniform} that there exists $L=L(\epsilon)$ such that uniform over large $p$ it holds that $\sum_{i>L} |a_{i,p}|^2 = O(\epsilon^2)$. Let
$$k_L(x,p) :=\sum_{h=0}^L a_{h,p} p_{h,p}(x)  \qquad \text{and} \qquad K_L(x,p) :=\sum_{h=0}^L a_{h,p} P_{h,p}(x) \ \ .$$ 
We obtain
\begin{eqnarray}
\label{e.kL}    \E [|k(\xi_p,p)-k_L(\xi_p,p)|^2] = O(\epsilon^2) \ \ , \\
\label{e.KL}   \E [|K(\xi_{G,p},p)-K_L(\xi_{G,p},p)|^2]  = O(\epsilon^2) \ \ .
\end{eqnarray}
Let $f_L$ and $F_L$ correspond to $k_L$ and $K_L$, and let $A_L$ and $G_L$ be generated from them respectively. 

It follows from \eqref{e.kL} and \eqref{e.KL} and Lemma~\ref{l.L2perturbation} below that
$$\E [m_A(z) -  m_{A_L}(z)] = O(\epsilon)   \ \ ,$$
$$\E [m_{G}(z) -  m_{G_L}(z)] = O(\epsilon)  \ \ .$$
The following Lemma is in \cite{cheng-singer}, to keep the paper self-contained, a short proof of Lemma~\ref{l.L2perturbation} is included.
\begin{lemma}\label{l.L2perturbation} Let $Y$ a random vector  of length $p$ whose coordinates are iid with mean $0$ and variance $1/p$. Let $Y'$ be an iid copy of $Y$. Assume that for $p$ large 
$$\E |f_1(Y^T Y',p) - f_2(Y^T Y',p)|^2 \le \epsilon^2/p \ \ . $$
Let $A_1$ and $A_2$  be generated from $f_1$ and $f_2$ using $n$ iid copies of $Y$. Then in the large $p$ large $n$ limit it holds that
$$\E |m_{A_1}(z) - m_{A_2}(z)| = O(\epsilon) \ \ .$$
\end{lemma}
\proof We largely follow \cite{cheng-singer}. Following the proof of \eqref{e.pointwise}, it is clear that
$$\E [|m_{A_1}(z)-m_{A_2}(z)|^2] \le C n^{-1} \sum_{i\ne j} \E [|A_{1,ij} - A_{2,ij}|^2]$$
$$= Cn^{-1} n(n-1) \E [|f_1(Y^T Y', p) - f_2(Y^T Y',p)|^2]$$
$$= O(\epsilon^2) \ \ ,$$
which implies the desired estimate.
\endproof
Thus, it suffices for \eqref{e.coupling2} to show that $\E [m_{A_L}(z)-m_{G_L}(z)] = O(\epsilon)$. This will be proved in two steps. First, it will be shown in Section~\ref{s.AtoG} that
\begin{equation}\label{e.coupling4}
\E [m_{A_L}(z)-m_{\widetilde G_L}(z)] = O(\epsilon) \ \ ,
\end{equation}             
where $\widetilde G_L$ is generated from $M_G$ using $k_L$ (as opposed to $K_L$, which was used to generate $G_L$). Then in Section~\ref{s.ktoK} it will be shown that  
\begin{equation}\label{e.kLtoKL}
\E [m_{G_L}(z)-m_{\widetilde G_L}(z)]  \to 0 \ \ .
\end{equation}

\subsection{Proof of \eqref{e.coupling4}: Conversion from $A_L$ to $\widetilde G_L$}\label{s.AtoG}

The proof of \eqref{e.coupling4} will follow the strategy in \cite{taovu12acta}; the idea is to  convert $A_L$ to $\widetilde G_L$ in $np$ steps, in each step one entry in $M_A$ is replaced by the corresponding entry in $M_G$. It suffices to show that in every step it holds that
$\E [\Delta m(z)] = O(n^{-2}\epsilon)$,
here $\Delta m$ is the difference between the Stieltjes transforms of the underlying matrices. This is the content of Lemma~\ref{l.inductstep}.

The generic setting for each step can be formulated as follows.  Let $M[1]$ and $M[2]$  denote two $p\times n$ random matrices that share the same entries except for the $(i,j)$ position. Assume that these entries are independent, and their distribution have mean $0$, variance $1/p$, and higher moments bounded (with uniform constants) by properly scaled powers of $p$. Let $A[1]$, $A[2]$ be generated from $M[1]$, $M[2]$ using the kernel function $f_L$. 

\begin{lemma} \label{l.inductstep}  In the large $n$ large $p$ limit
$$\E  [m_{A[1]}(z) -  m_{A[2]}(z)]    = O(n^{-2}\epsilon) \ \ .$$
\end{lemma}
For simplicity of notation, in the rest of this section $M[0]$ denotes the matrix that shares the same entries with $M[1]$, $M[2]$ except for the $(i,j)$ position, where  $M[0]_{ij}:=0$. Denote by $A[0]$ the kernel matrix generated from $M[0]$ using $f_L$. 

\proof   Let $E$ be the event that $\|A[0]-A[m]\| \le Im(z)/2$ for both $m=1,2$. In Lemma~\ref{l.momentbound}, it will be shown that 
\begin{equation}\label{e.momentbound}
\E [ \|A[0]-A[m]\|^{q}] = O(n^{-q/2})  \ \ , \ \  m = 1 , 2 \ \ ,
\end{equation}  
here $\|.\|$ denotes the matrix norm. Since $q$ could be taken large ($q>4$ suffices), it follows that $E^c$  has probability $o(n^{-2})$. 
Clearly $m_{A[m]}(z)=O(1)$, thus it remains to show
\begin{equation}\label{e.inductstep}
\E \Big(1_{E } [m_{A[1]}(z) -  m_{A[2]}(z)]\Big) = O(n^{-2}\epsilon)
\end{equation}
in the large $n$ large $p$ limit.

Let $A$ denote $A[1]$ or $A[2]$ in the rest of the proof, and let $M \in \{M[1], M[2]\}$ be the corresponding sample matrix. On $E$, expand
\begin{equation}\label{e.expansion}
(A-zI)^{-1} = \sum_{k\ge 0} R_k \ \  ,  \ \  \text{where}
\end{equation}
$$R_k:=[(A[0]-zI)^{-1}(A[0]-A)]^k (A[0]-zI)^{-1} \ \ .$$
In particular, $R_0=(A[0]-zI)^{-1}$. Note that any eigenvalue of $A[0]-Iz$ is of the form $\lambda-z$ for some eigenvalue $\lambda$ of $A[0]$. Since any such $\lambda$ is real, it follows that
\begin{equation}\label{e.resolvent}
\|R_0\| \le \frac1 {Im(z)} \ \  ,
\end{equation}
Thus on $E$ the expansion \eqref{e.expansion} is absolutely convergence with respect to $\|.\|$.

This expansion will be used to compute the leading `asymptotics' of $\E [1_{E} m_A(z)]$ as $n\to\infty$, and show that modulo $o(n^{-2})$ one obtains the same leading asymptotics for   $A=A[1]$ and $A=A[2]$, and this clearly implies \eqref{e.inductstep}.

In \eqref{e.expansion}, it is clear that the contributions to $\E [1_{E} m_A(z)]$ of $R_0$ is the same for $A=A[1]$ and $A=A[2]$, so below only $R_k$ with $k\ge 1$ is considered. \\

\noindent \underline{\emph{Decay estimates for contribution of higher order terms in \eqref{e.expansion}}}:\\
First, it will be shown that  for some absolute constant $C$ the following holds for $k\ge 1$ :
\begin{equation}\label{e.tail}
|\Tr (R_k)|  \le C Im(z)^{-(k+1)} \|A[0]-A\|^k
\end{equation}
To show \eqref{e.tail}, the key observation is that $A[0]-A$ has rank at most $2$. Indeed, this follows from the fact that at most one column and at most one row in $A[0]-A$ could be nonzero. Thus, $[R_0(A[0]-A)]^k R_0$ is of rank at most $2$ and therefore
$$|\Tr (R_k)|  \le 2 \Big\|[R_0(A[0]-A)]^k R_0 \Big\|$$
\begin{equation}\label{e.tailresolvent}\le 2 \|R_0\|^{k+1} \|A[0]-A\|^{k} \ \ ,
\end{equation}
and \eqref{e.tail} follows immediately from \eqref{e.tailresolvent} and \eqref{e.resolvent}. 

As a consequence of \eqref{e.tail}, it follows that
$$n^{-1}\E \Big[1_{E} \,  |\sum_{k\ge 3}\Tr R_k| \Big] \le C n^{-1} \, \E  [ \|A[0]-A\|^{3}] = O(n^{-5/2}) \ \ .$$

\noindent \underline{\emph{Asymptotics matching for the contribution of $R_1$}}:\\
Rewrite
$$n^{-1}\E[1_{E}\Tr (R_1)] $$
\begin{equation}\label{e.Edecomp}
= n^{-1}\E[ \Tr (R_1)] - n^{-1} \E[1_{E^c}\Tr (R_1)] \ \ .
\end{equation}
We first show that the second term in \eqref{e.Edecomp} is $o(n^{-2})$. Indeed, it follows from Cauchy-Schwarz that
$$n^{-1}|\E [1_{E^c} \Tr (R_1)]|$$
$$ \le n^{-1}\P(E^c)^{1/2} (\E [|\Tr R_1|^2])	^{1/2}$$
$$\le Cn^{-1}\P(E^c)^{1/2} (\E [\|A[0]-A\|^2] )^{1/2} \qquad \text{(using \eqref{e.tailresolvent})}$$
$$\le C n^{-1}o(n^{-1}) O(n^{-1/2}) = o(n^{-5/2})$$
which implies the desired estimate.

It remains to compute the asymptotic of the first term in \eqref{e.Edecomp}. Let $\E_0$ denote expectation with respect to entries of $M[0]$ and let $\E^{(ij)}$ denote expectation with respect to the $(i,j)$ entries of $M[1]$ and $M[2]$.  Since $A[0]$ is independent of the $(i,j)$ entries, it follows that
$$\E[ \Tr (R_1)] = \Tr \,  \, \E_0 \Big[R_0 \,\Big\{ \E^{(ij)} \Big([A[0]-A]\Big) \Big\} \, R_0 \Big]  \ \ . $$
It will be shown in Lemma~\ref{l.entrydecomposition} that there exists a decomposition
\begin{equation}\label{e.entrydecomposition}
A[0] - A = a M_{ij}    + b M_{ij}^2   + c  \ \ ,
\end{equation}
where $a,b,c$ are $n\times n$ matrices such that the following holds:
\begin{itemize}
\item The entries outside the $j^{th}$ rows and the $j^{th}$ columns of $a$, $b$, $c$ are zeros.
\item The nonzero entries of $a$ and $b$ depend only on $f$ and $M[0]$, 
\item the $2^{nd}$ moment of any entries of $a$ and $b$ are of size $O(n^{-1})$,  and  the $2^{nd}$ moment of any entries of $c$ are of size $O(n^{-4})$.
\end{itemize}
It follows that $\E[\|c\|^2] = O(n^{-3})$ and $c$ has rank $2$, hence
$$\Tr [\E_0 R_0 (\E^{(ij)} c) R_0 \Big]  = \Tr \E \Big[R_0 c R_0 \Big] = O(\|c\|) = O(n^{-3/2}) \ \ .$$
It follows that
$$n^{-1}\E[\Tr (R_1)]  = n^{-1} \E_0\Big(\Tr\, \big[R_0 a R_0\big]\Big)  \E^{(ij)} [M_{ij}] $$
$$+ n^{-1}\E_0\Big(\Tr\, \big[R_0 b R_0\big]\Big)  \E^{(ij)} [M_{ij}^2] $$
$$+ O(n^{-5/2}) \ \ .$$
We remark that the first two terms are the same for $A=A[1]$ and $A=A[2]$.  \\

\noindent \underline{\emph{Asymptotics matching for the contribution of $R_2$}}:\\
As before, rewrite
$$n^{-1} \E[1_E \Tr(R_2)] = n^{-1} \E[\Tr(R_2)] - n^{-1} \E[1_{E^c} \Tr(R_2)]$$
and using \eqref{e.tailresolvent} it is not hard to see that the second term is $o(n^{-2})$. Therefore it remains to compute the asymptotics for the first term. Again, the decomposition \eqref{e.entrydecomposition} will be used to expand
\begin{equation}\label{e.R2expansion}
n^{-1} \E [\Tr (R_2)] = n^{-1} \E^{(ij)} [M_{ij}^2] \, \E_0 \Big( \Tr[(R_0 a)^2 R_0] \Big)+ \text{other terms} \ \ .
\end{equation}
Since the second moment of any entries of $a$ and $b$ are $O(n^{-1})$ and since $a$ and $b$ have $O(n)$ non-zero entries, it is not hard to see that $\E[\|a\|^2] = O(1)$ and $\E[\|b\|^2]= O(1)$, while $\E[\|c\|^2] = O(n^{-3})$ as remarked above. Therefore using the small rank properties of $a$, $b$, $c$, the \emph{other terms} in  the expansion \eqref{e.R2expansion} can be bounded  by expected values of products of spectral norms, and eventually obtain an estimate of $O(n^{-5/2})$. Clearly, the first term in the expansion \eqref{e.R2expansion} is the same for $A=A[1]$ and $A=A[2]$.

Finally, to complete the proof of Lemma~\ref{l.inductstep}, it remains to show \eqref{e.momentbound} (see Lemma~\ref{l.momentbound}) and \eqref{e.entrydecomposition} (see Lemma~\ref{l.entrydecomposition}).
\endproof

To keep the following Lemmas self-contained, the symbol $H[k]$ will be used instead of $A[k]$. In the applications of these Lemma to obtain \eqref{e.momentbound}  and  \eqref{e.entrydecomposition}, these are the same.

\begin{lemma}\label{l.momentbound} For each $0\le k \le 2$ let $H[k]$ be generated from the column vectors of $M[k]$ using a kernel function $h(x,p)$ using the non-diagonal model \eqref{e.Adefnodiag}. Assume that
\begin{equation}\label{e.h_x}
|\partial_x h(x,p)| \le C_N (1+|\sqrt p x|^N) \ \ ,
\end{equation}
for some $N>0$ uniform over $p$ and $x$, Then for $m=1,2$ and $q>0$ it holds that
$$\E \Big(\|H[0]-H[m]\|^{q}\Big) = O_N (p^{-q/2}) \ \ .$$
\end{lemma}
Remark: To obtain \eqref{e.momentbound}, Lemma~\ref{l.momentbound} will be applied for $h(x,p)=f_L(x,p)$, which satisfies \eqref{e.h_x} thanks to Lemma~\ref{l.opcompare} and Lemma~\ref{l.opcoefficients}. Note that in that case it holds that $H[k]=A[k]$.

\proof For simplicity, let $H=H[m]$ and $M=M[m]$ and  $U=M[0]^T M[0]$.  Since $H[0]-H$ has at most $O(n)$ nonzero entries and since 
$$\|B\| \le (\sum_{i,j} |B_{ij}|^2)^{1/2}$$
for any square matrix $B$, it suffices to show that the $q$ moment of each nonzero entry of $H[0]-H$ is bounded above by $O(n^{-q})$.  

Consider without loss of generality an (non-diagonal) entry on the $j^{th}$ row of $H[0]-H$. This entry has the form
$$\alpha_{jk}=h(U_{jk},p) - h(U_{jk} + M_{ij} M_{ik},p)$$
for some $1\le k \le n$. It follows from \eqref{e.h_x} that
$$|\alpha_{jk}|  = O_N\Big(|M_{ij}M_{ik}| \big(1+\sqrt p |M_{ij}M_{ik}| + \sqrt p |U_{jk}| \big)^N\Big)$$
Using \eqref{e.momentassumption}, it is not hard to see that $E \big[|\sqrt p U_{jk}|^K \big]= O(1)$ (for details see for instance the proof of \eqref{e.innerproduct}). Also, it is clear that $\E |pM_{ij} M_{ik}|^K = O(1)$. Therefore, using Cauchy-Schwarz it follows that
$$\E \Big[|\alpha_{jk}|^q \Big] \le C_N \Big(\E |M_{ij}M_{ik}|^{2q}\Big)^{1/2} = O(p^{-q}) \ \ .$$
\endproof

\begin{lemma}\label{l.entrydecomposition} For each $0\le k \le 2$ let $H[k]$ be generated from the column vectors of $M[k]$ using a kernel function $h(x,p)$ using the non-diagonal model \eqref{e.Adefnodiag}.
Assume that for each $m=1,2,3$ there is some $N>0$  such that
\begin{equation}\label{e.h_xxx}
|\frac {\partial^m}{\partial x^m} h(x,p)| \le C_N  p^{(m-1)/2}(1+|\sqrt p x|^N) \ \ ,
\end{equation}
where $C_N$ is uniform over $p$ and $x$. Then for $m=1,2$ it holds that
$$ H[m]  - H[0]  = a M[m]_{ij}    + b \Big[M[m]_{ij}\Big]^2   + c  \ \ ,$$
where $a$, $b$, $c$ are $n\times n$ matrices such that 
\begin{itemize}
\item[(i)] the entries outside the $j^{th}$ rows and columns of $a$, $b$, $c$ are zeros, and
\item[(ii)] the nonzero entries of $a$ and $b$ are independent of $M[m]_{ij}$, and
\item[(iii)] the $2^{nd}$ moment of any entries of $a$ and $b$ are of size $O(n^{-1})$,  and  the $2^{nd}$ moment of any entries of $c$ are of size $O(n^{-4})$.
\end{itemize}
\end{lemma}
Remark: The condition \eqref{e.h_xxx} is satisfied when $h=f_L$ thanks to Lemma~\ref{l.opcompare} and Lemma~\ref{l.opcoefficients}. Thus, \eqref{e.entrydecomposition} follows by applying Lemma~\ref{l.entrydecomposition} to $h=f_L$, in that case it holds that $H[k]=A[k]$.

\proof Let $H=H[m]$ and $M=M[m]$ for $m\in \{1,2\}$. Let $U=M[0]^T M[0]$. It is clear that all nonzero entries of $H-H[0]$  must be in the $j^{th}$ column or the $j^{th}$ row and off the diagonal, thus it remains to decompose these non-zero entries. 

Consider without loss of generality an entry in the $j^{th}$ row of $H-H[0]$, which has the form
$$\alpha_{jk}=  h(U_{jk} + M_{ij} M[0]_{ik},p) - h(U_{jk},p)$$
for some $k\ne j$. We now decompose
$$\alpha_{jk} =  a_{jk} M_{ij}  + b_{jk} M_{ij}^2 + c_{jk}$$
$$a_{jk} =   M[0]_{ik} h_x(U_{jk},p)  \ \ , \ \ b_{jk} =   \frac 1  2 \Big[M[0]_{ik} \Big]^2 h_{xx}(U_{jk},p) \ \ .$$
Similar decomposition for the $j^{th}$ column of $A[0]-A$ also holds. It is clear that the two conditions (i) and (ii) are satisfied and it remains to show that
$$\E |c_{jk}|^2 = O(n^{-4}) \ \ .$$
By given assumption, for some $\theta$ in between $U_{jk}$ and $U_{jk} + M_{ij}M[0]_{ik}$ it holds that
$$ c_{jk}  = \frac 1 6 \Big[M_{ij} M[0]_{ik}\Big]^3  h_{xxx}(\theta,p)  \ \ .$$
In particular, it follows from \eqref{e.h_xxx} that, for some $N>0$ and $C$ depends on $N$,
$$|c_{jk}| \le C \, p |M_{ij} M[0]_{ik}|^3 \Big[1+ (\sqrt p |\theta|)^N\Big]  \ \ ,$$
$$\le C\,  p |M_{ij} M[0]_{ik}|^3 \Big[1+ p^{N/2}\big(|U_{jk}| + |M_{ij} M[0]_{ik}|\big)^N\Big]   \ \ .$$
Using H\"older's inequality and the previously obtained bounds on moments of entries of $M$, $M[0]$, $U$, it follows that
$$\E |c_{jk}|^2 = O(p^{-4}) \ \ .$$
\endproof

\subsection{Proof of \eqref{e.kLtoKL}: Conversion from $\widetilde G_L$ to $G_L$}\label{s.ktoK}
Recall that $\xi_{G,p}=\sqrt p Y^T Y'$ where $Y$ and $Y'$ denote  iid $1\times p$ vectors with independent Gaussian $N(0,1/p)$ coordinates. Applying Lemma~\ref{l.L2perturbation}, it remains to show that: for any $\delta>0$, the following holds uniformly over $p$ large
\begin{equation}\label{e.kLtoKLinL2}
\E \big[ |k_L(\xi_{G,p}, p) - K_L(\xi_{G,p},p)|^2 \big] = O(\delta) \ \ ,
\end{equation}
here the implicit constant depends on $L$. 

It follows from the triangle inequality and Lemma~\ref{l.opcompare} that, for $p$ large, 
$$|p_{k,p}(x) - P_{k,p}(x)|  \le C_k\delta(1+|x|^k) \ \ .$$
On the other hand it is clear that  $\E |\xi_p|^k = O_k(1)$. It follows that, for $p$ large, 
$$\|p_{k,p}-P_{k,p}\|_{L^2(\mu_{G,p},\R)} \le C_k \delta $$
for each $k=\overline{0,L}$.
Now \eqref{e.kLtoKLinL2} follows since 
$$K_L(x,p)  - k_L(x,p) =  \sum_{k=1}^L a_{k,p} \Big[P_{k,p}(x)-p_{k,p}(x)\Big] \ \ ,$$
and  for $p$ large  $a_{k,p} = O(1)$ by Lemma~\ref{l.opcoefficients}.

\section{Concluding remarks}

In Theorem~\ref{t.main}, it was assumed that the entries of $X_i$'s are i.i.d., while  in Theorem~\ref{t.maincor1} and Theorem~\ref{t.maincor2} it is possible to have random vectors with dependent entries, as long as a high concentration condition is satisfied. Cheng and Singer \cite{cheng-singer} on the other hand have outlined a proof of an analogue of Theorem~\ref{t.main} in the setting when $X_i$'s are independently identically sampled from the unit sphere. This suggests that the i.i.d. assumption on entries of $X_i$'s may be weakened. Our proof of Theorem~\ref{t.main} in this paper however relies on the independence of entries of $X_i$, more specifically in the implementation of the Lindeberg swapping argument in Section~\ref{s.AtoG}.  It would be interesting to see if this swapping argument could be improved to extend Theorem~\ref{t.main} to settings when $X_i$'s have dependent entries. In this direction, see for instance Chatterjee \cite{chatterjee} where some generalization of the Lindeberg principle was considered.

In a different direction, one may ask questions about local statistics of the eigenvalues of random kernel matrices. The local statistics of Wigner and covariance matrices have been studied extensively in the literature, see e.g. \cite{taovu12acta} or \cite{taovusurvey} for a comprehensive survey. However, we are not aware of any related work in the setting of random kernel matrices, even when the envelop function is independent of $p$. It seems that a naive adaptation of the approximation argument, carried out in this paper and El Karoui's work \cite{karoui07a, karoui07b, karoui}, does not lead to sufficiently interesting information about local statistics of the eigenvalues, unless very special assumptions are made on the kernel. El Karoui \cite{karoui07a, karoui07b, karoui} on the other hand has been able to obtain some results about behavior of the largest eigenvalues via the approximation approach.

\end{document}